\documentclass[11pt]{article}

\usepackage{amsmath,amssymb}
\usepackage{latexsym}
\usepackage{graphicx}
\usepackage{bm}

\newtheorem{theorem}{Theorem}

\begin{document}

\markboth{N. S. Hoang and A. G. Ramm}{On stable numerical differentiation}

\title{
On stable numerical differentiation}

\author{N. S. Hoang$\dag$ and A. G. Ramm$\dag$\footnotemark[1]\\
\\
$\dag$Mathematics Department, Kansas State University,\\
Manhattan, KS 66506, USA
}

\date{}
\maketitle

\renewcommand{\thefootnote}{\fnsymbol{footnote}} 
\footnotetext[1]{
Corresponding author.
\\email: nguyenhs@math.ksu.edu (N. S. Hoang), \, ramm@math.ksu.edu (A. G. Ramm)
} 

\begin{abstract}
\noindent
Based on a regularized Volterra equation,
two different approaches for numerical differentiation are considered. 
The first approach consists of solving a regularized Volterra equation 
while
the second approach is based  on solving a disretized version of the 
regularized 
Volterra equation.
Numerical experiments show that these methods are efficient and 
compete favorably
with the variational regularization method for stable
calculating the derivatives of noisy functions. 


\end{abstract}

\section{Introduction\label{s1}}

Calculating the derivatives of noisy functions is of prime importance
in many applications. 
The problem consists of calculating stably  the derivative of a smooth 
function $f$
given its noisy data $f_\delta$, $\|f_\delta-f\|\leq \delta$.
This is an ill-posed problem: a small error in $f$ may lead to
a large error in $f'$.
Many methods have been introduced in the literature. A review is given in [7].
Divided differences method with $h=h(\delta)$ has been first proposed in [4],
see also \cite{133, 168, 181}.
Necessary and sufficient conditions for the existence of a method
for stable differentiation of noisy data are given in 
\cite[chapter 15]{Ramm}, see also \cite{441}.
In our paper a method for stable differentiation
based on solving the regularized Volterra equation
\begin{equation}
\label{eq1}
Au(x)+ f_\delta(0):=\int_0^x u(s)ds + f_\delta(0)=f_\delta(x),
\end{equation}
is proposed (see also \cite{415,494,441}). 
One often
applies the Variational Regularization (VR) method
\begin{equation}
\label{eqregu}
\|Au-f_\delta\|^2 + \alpha \|u\|^2 \to \min
\end{equation}
for stable differentiation.

 In this paper (and in \cite{494}) an approach, based 
on the fact that the quadratic form 
of the operator $A$ is nonnegative in real Hilbert space $L^2(0,a)$, 
$a=const>0$, is used.

\section{Methods\label{s2}}

Consider two different approaches to solving equation \eqref{eq1}. The first 
approach consists of
solving  directly regularized equation \eqref{eq1}.
The second approach is based on the Dynamical Systems method (DSM) and 
an iterative scheme  from \cite{526}.

\subsection{First method}

In \cite{494}, the derivatives of a noisy function $f_\delta$ 
are obtained by solving the equation
\begin{equation}
\label{eq2}
\alpha u_{\alpha,\delta} + A u_{\alpha,\delta} = f_\delta.
\end{equation}
If $\alpha=\alpha(\delta)>0$ is continuous on $[0,\delta_0)$, 
$\delta_0>0$ and  
\begin{equation}
\label{eq3}
\lim_{\delta\to 0}\alpha(\delta)=0,\quad \lim_{\delta\to0}\frac{\delta}{\alpha(\delta)}=0,
\end{equation}
then the following result holds (see \cite{494}): 

\begin{theorem}
Assume \eqref{eq3}. Then
$$
\lim_{\delta\to0}\|u_\delta - u\|=0,
$$
where $u_\delta$ solves \eqref{eq2} with $\alpha=\alpha(\delta)$.
\end{theorem}

The solution of \eqref{eq2} is:
\begin{equation}
\label{eq4}
u_\delta(x) = -\frac{1}{\alpha^2}\exp(-\frac{x}{\alpha})
\int_0^x\exp(\frac{s}{\alpha})f_\delta(s)ds + \frac{f_\delta(x)}{\alpha}.
\end{equation}
This formula  and an {\it a priori} choice $\alpha(\delta) = \delta^k/c$, 
where $k\in (0,1)$, $c$ is a constant,  
yield a scheme for stable differentiation.
When $\alpha(\delta)$ is known, 
the problem is reduced to 
calculating integral \eqref{eq4}.
There are many methods for calculating
accurately and fast integral \eqref{eq4} (see e.g. \cite{Davis}).
However, there is no known algorithm for choosing $k,c$ which are optimal in some sense.
The advantage of our approach is that the CPU time for the method 
is very small compared with the VR and DSM, see  Section \ref{sectionfirst}. 
Moreover, one can calculate the solution
analytically when the function $f_\delta$ is simple by using  tables of integrals or 
MAPLE.

\subsection{An iterative scheme of DSM for solving discretizations of
the regularized Volterra equation}
\label{iterativesec}

Another approach to stable differentiation 
is to use the DSM (see \cite{Ramm}).
The DSM yields a stable solution of the equation:
\begin{equation}
\label{eq5}
F(u)=Au-f=0,\quad u\in H,
\end{equation}
where $H$ is a Hilbert space and $A$ is a linear operator in $H$ which is 
not necessarily bounded but closed and densely defined.
The DSM to solve \eqref{eq5} is of the form:
\begin{equation}
\label{eq6}
u'=-u+(T+a(s))^{-1}A^{*}f,\quad u(0)=u_0,
\end{equation}
where $T:=A^{*}A$ and $a(t)>0$ is a nonincreasing function such that 
$a(t)\to 0$ as $t\to\infty$. 
The unique solution to \eqref{eq6} is given by
\begin{equation}
\label{eq7}
u(t)=u_0e^{-t}+e^{-t}\int_0^te^{s}(T+a(s))^{-1}A^{*}fds.
\end{equation}
An iterative scheme
for computing $u(t)$ in \eqref{eq7} is proposed in \cite{526}:
$$
u_{n+1}=e^{-h_n}u_{n}+(1-e^{-h_n})\big{(}T+a_n\big{)}^{-1}A^*f_\delta,\quad h_n=t_{n+1}-t_n.
$$
With $a_0$ satisfying 
\begin{equation}
\label{a0cond}
\delta<\|Au_{a_0} - f_\delta\| < 2\delta,
\end{equation}
one chooses $a_n$ and $h_n$ as follows:
$$
a_n = \frac{a_0}{1 + t_n},\quad h_n = q^n,
$$ 
where $1\leq q\leq 2$, $t_0 = 0$.
To increase the speed of computing we recommend choosing $q=2$.
At each iteration one checks if
\begin{equation}
\label{eq8}
0.9\delta \leq \|Au_n-f_\delta\|\leq 1.001 \delta.
\end{equation}
This is a stopping criterion of discrepancy principle type (see 
\cite{526}).
If $t_n$ is the first time such that \eqref{eq8} is satisfied, then one stops and
takes $u_n$ as the solution to \eqref{eq5}. 
The choice of $a_0$ satisfying \eqref{a0cond} is done by iterations as follows:

\begin{enumerate}
\item{As an initial guess for $a_0$ one takes 
$a_0=\frac{1}{3}\|A\|^2\delta_{rel}$, where 
$\delta_{rel}=\frac{\delta}{\|f\|}$.} 
\item{
If  $\frac{\|Au_{a_0} -f_\delta\|}{\delta}=c>3$, 
then one takes $a_1:=\frac{a_0}{2(c-1)}$ as the next guess and checks 
if condition \eqref{eq8} is satisfied. 
If $2<c\leq 3$ then one takes $a_1:=a_0/3$.}
\item{
If  $\frac{\|Au_{a_0} -f_\delta\|}{\delta}=c<1$, then $a_1:=3a_0$ is used as the next guess.}
\item{
After $a_0$ is updated, one checks if \eqref{eq8} is satisfied. 
If \eqref{eq8} is not satisfied, one repeats steps 2 and 3 until 
one finds $a_0$ satisfying condition \eqref{eq8}.} 
\end{enumerate}
Algorithms for choosing $a_0$ and computing $u_n$ are detailed 
in algorithms 1 and 2 in \cite{526}.

\section{Numerical experiments\label{s3}}

Numerical experiments are carried out in MATLAB in double-precision
arithmetic. In all experiments, by $u(t)$,
$u_{[1]}(t)$, $u_{\text{DSM}}(t)$ and $u_{\text{VR}}(t)$ we denote the 
exact derivative, the derivatives computed by the first, the DSM
and the VR methods, respectively. In this section by $n$ we denote the number of points 
used to discretize the interval $[0,1]$.

\subsection{Computing the first derivatives of a noisy function}
\label{sectionfirst}

\begin{figure}[!h!tb]
\centerline{%
\includegraphics[scale=0.85]{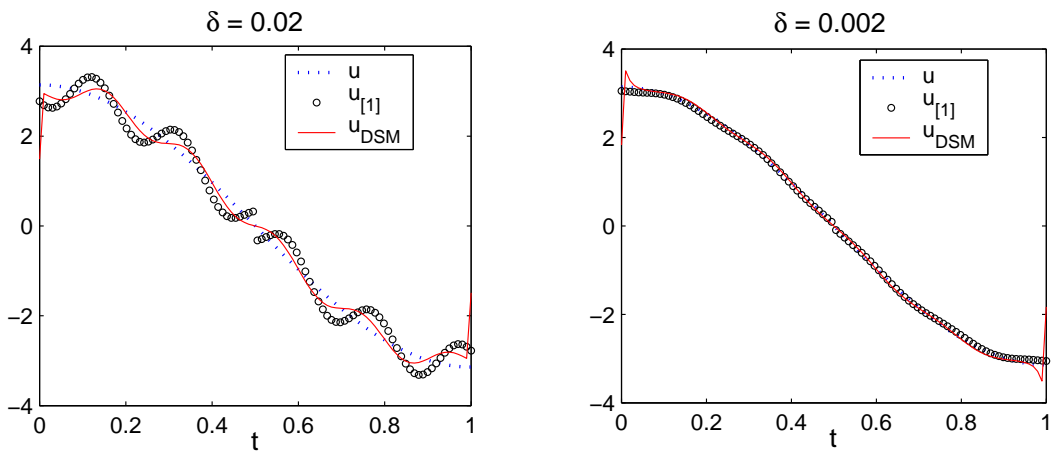}}
\caption{Numerical results for $f_\delta(x)=\sin(\pi t)+\delta \cos(10\pi t)$. 
Discretization points $n=100$.}
\label{fig1}
\end{figure}

\begin{figure}[!h!tb]
\centerline{%
\includegraphics[scale=0.85]{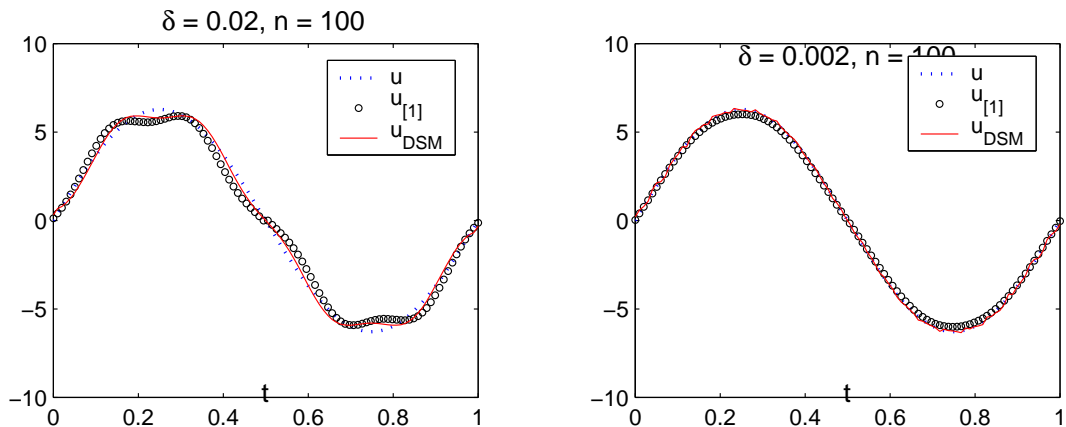}}
\caption{Numerical results for $f_\delta(x)=\sin(2\pi t - \frac{1}{2}\pi)+\delta \cos(10\pi t)$. 
Discretization points $n=100$.}
\label{fig12}
\end{figure}

Let us compute the derivatives of the function $f(t) = \sin(\pi t)$ 
contaminated by the noise function $e(t)= \delta \cos(10\pi t)$. 
The derivative of $f(t)$ is $f'(t)=\pi \cos(\pi t)$. 
To solve this problem we use three methods: the first method, 
based on computing
integral \eqref{eq4}, the VR method, 
and the DSM method, based on a discretized version of 
\eqref{eq1}. Numerical 
results for this problem are presented in Figure~\ref{fig1}. In our experiments,
since the results otained by the DSM and the VR are nearly the same, we 
present only the results for the DSM in Figure~\ref{fig1} and \ref{fig12} in order to make 
these figures simple.

In this experiment the trapezoidal quadrature rule is applied to integral 
equation \eqref{eq1} and is used for computing integral \eqref{eq4}. 
One may use higher order intepolation methods to compute integral \eqref{eq4}.
However, it does not necessarily bring improvements in accuracy. 
This is so because using 
a high order intepolation method for inaccurate data may even lead to worse results.
This is the case when the noise level is large. 

The approximate derivative formula \eqref{eq4} for $t$ close to 0
does not use much information about $f_\delta$. Thus, we only use \eqref{eq4}
for computing $f'(t)$ for $t\in [\frac{1}{2},1]$. For $t\in 
[0,\frac{1}{2})]$,
we take $g_\delta(t):=f_\delta(1-t)$ and use formula \eqref{eq4} for
 $g_\delta(t)$
with $t\in (\frac{1}{2},1]$. That is, we have a discontinuity at 
$t=\frac{1}{2}$ of
the solution, obtained by the first method in Figure~\ref{fig1} 
and \ref{fig12}.
The same idea is applied in discretizing equation \eqref{eqregu} in 
the implementation
of the DSM and VR.

In the DSM and VR we also use the trapezoidal quadrature rule to 
discretize equation \eqref{eq1}.
Since the right-hand side $f_\delta$ contains noise, using high order 
collocation methods
does not necessarily improve the accuracy. 
Experiments have shown that the use of higher order collocation methods 
leads to linear algebraic systems with larger condition numbers and 
yields  
numerical solutions with low accuracy. 

The CPU times for the VR  and DSM are about 0.0125 sec.
The  CPU time for the first method is much smaller: 0.0015 sec.
Here, we should bear in mind that the DSM and the VR use iterations
to look for "good" regularization parameter $\alpha$ while the code based on the first method
does nothing to look for $\alpha$ but uses $\alpha$ as an input value.
If one also uses the regularization parameter as an input in the VR and DSM, although these
methods still take more time than the first method the difference in
computation time is not so large. 

The error of the first method  for $\delta=0.02$ is larger than 
those of the VR and the DSM, but when $\delta=0.002$ then 
the first method gives smaller errors. From Figure~\ref{fig1} and \ref{fig12},
one can see that the solutions obtained by the DSM are better 
than those obtained by the first method
for all $t\in [0,1]$ except for the $t$ which are close to the boundary of 
the interval.
Indeed, it can be showed analytically that the solution $u$ to equation \eqref{eqregu} satisfies $u(0)=u(1)=0$.
However, the derivative of $f$ in Figure~\ref{fig1} satifies $f'(0)=\pi$ and $f'(1)=-\pi$. 
If the computed derivatives at the points close to the boundary 
are discarded,  then in both cases the DSM and the VR are more accurate 
than the first method.

Figure~\ref{fig12} presents the numerical experiment
for $f(t)=\sin(2\pi t - \frac{1}{2}\pi)$ contaminated by the same noise function $e(t)= \delta \cos(10\pi t)$.
For this problem, since
the function to be differentiated $f$ satisfies $f'(0)=f'(1)=0$ both the 
DSM and the VR give more accurate
results than the first method.

From Figure~\ref{fig1} and ~\ref{fig12} one can see
that for $\delta = 0.02$ the computed derivatives are very close to the
exact derivative at all points except for those close to the boundary in Figure~\ref{fig1}.

\subsection{Computing the second derivatives of a noisy function}

Let us give numerical results for computing the second derivatives of noisy functions.
The problem is reduced to an integral equation of the first kind.
A linear algebraic system is obtained by a discretization of the 
integral equation whose kernel $K$ is Green's function
$$
K(s,t)=
\left\{
\begin{matrix}
s(t-1),\quad \text{if}\quad s<t\\
t(s-1),\quad \text{if}\quad s\geq t
\end{matrix}
\right. .
$$
Here $s,t\in[0,1]$ and as the right-hand side $f$ and the corresponding 
solution $u$ one chooses one of the following 
(see \cite{526}):
\begin{align*}
\text{case 1},\quad &f(s)=\frac{s^3-s}{6},\quad u(s)=s,\quad 0\leq s\leq 1,\\
\text{case 2},\quad &f(s)=\frac{\sin(2\pi s)}{4\pi^2}+s-1,\quad u(s)=\sin(2\pi s),\quad 0\leq s\leq 1.
\end{align*}

Collocation method is used for discretization. 
This discretization can be improved by other methods 
but we do not go into detail.
We use $n=10,20,...,100,$ and $b_{n,\delta} = b_n + e_n$, 
where $e_n$ is a vector containing random entries, normally distributed 
with mean 0, variance 1, and scaled so that $\|e_n\|=\delta_{rel}\|b_n\|$. 
This linear algebraic system is mildly ill-posed: the condition 
number of $A_{100}$ is $1.2158\times10^{4}$. 

\begin{table}[h] 
\caption{Results for case 1 and 2 with $\delta_{rel}=0.01$, $n=20,40,...,100$.}
\label{table1}
\centering
\small
\begin{tabular}{@{  }c@{\hspace{2mm}}|c@{\hspace{2mm}}c@{\hspace{2mm}}
|c@{\hspace{2mm}}c@{\hspace{2mm}}|c@{\hspace{2mm}}|c@{\hspace{2mm}}
c@{\hspace{2mm}}|c@{\hspace{2mm}}c@{\hspace{2mm}}
|c@{\hspace{2mm}}c@{\hspace{2mm}}|c@{\hspace{2mm}}r@{\hspace{2mm}}l@{}
} 
\hline
\multicolumn{5}{c|}{Case 1}&\multicolumn{5}{c|}{Case 2}\\
\hline
&\multicolumn{2}{c|}{DSM}&\multicolumn{2}{c|}{VR}& &\multicolumn{2}{c|}{DSM}&\multicolumn{2}{c|}{VR}\\
$n$&
$N_{\text{linsol}}$&$\frac{\|u_\delta-y\|_{2}}{\|y\|_2}$&
$N_{\text{linsol}}$&$\frac{\|u_\delta-y\|_{2}}{\|y\|_2}$&
$n$&
$N_{\text{linsol}}$&$\frac{\|u_\delta-y\|_{2}}{\|y\|_2}$&
$N_{\text{linsol}}$&$\frac{\|u_\delta-y\|_{2}}{\|y\|_2}$
\\
\hline
 20    &3    &0.3319    &5    &0.3440	&20	&4    &0.0773	&4    &0.0780\\
 40    &4    &0.3206    &6    &0.3253	&40	&3    &0.0484	&6    &0.0520\\
 60    &4    &0.3264    &6    &0.3312	&60	&4    &0.0355	&6    &0.0438\\
 80    &4    &0.3019    &7    &0.3014	&80	&3    &0.0407	&5    &0.0479\\
100    &5    &0.2956    &7    &0.2948	&100	&4    &0.0254	&6    &0.0379\\
\hline 
\end{tabular}
\end{table}

\begin{figure}[!h!tb]
\centerline{%
\includegraphics[scale=0.85]{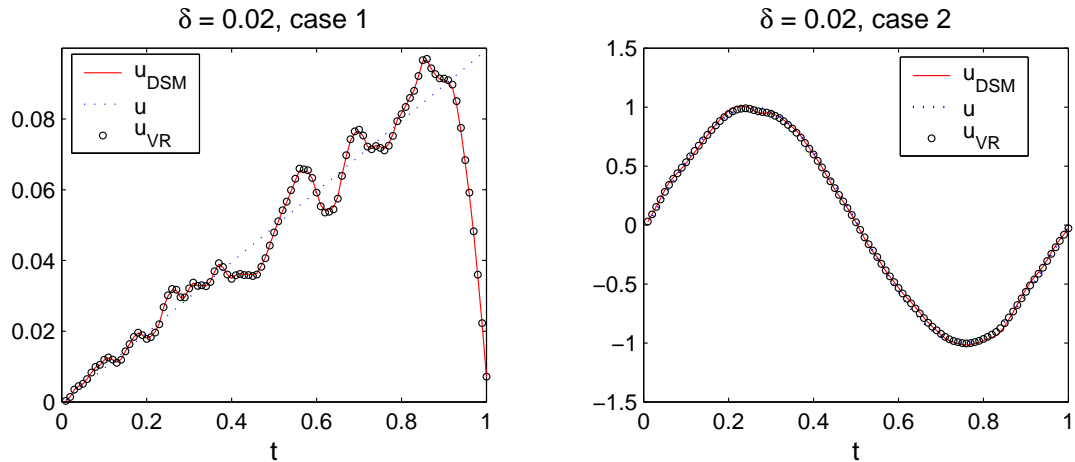}}
\caption{Plots of solutions obtained by DSM, VR when $n=100$, 
$\delta_{rel}=0.02$.}
\label{fig2}
\end{figure}

Table~\ref{table1} shows that
numerical results obtained by the DSM are  
more accurate than those by the VR. 
Figure~\ref{fig2}  plots the numerical solutions for these cases. 
The computation time of the DSM in these cases is about the same as 
or less than that of the VR.
From Table~\ref{table1} 
one can see that both the DSM and the VR perform better 
in case 2 than in  case 1.
Note that the regularized equation to solve for second derivatives in this case is
of the same form as equation \eqref{eqregu}.
As we discussed earlier, it is because in case 2 we have $f'(0)=f'(1)=0$.

We conclude that in this experiment the DSM competes favorably 
with the VR.

Looking at Figure~\ref{fig2} case 1, one can see that the computed 
values at $t=0$ and $t=1$ are zeros.
Again, the regularized scheme forces the computed derivative $u$ to 
satisfy the relations $u(1)=u(0)=0$. 
If one wants to compute the derivative of a noisy function on an 
interval by the proposed method, one should collect data on a larger 
interval and use this method to calculate the derivative 
at the points which are not close to the boundary.

\section{Concluding remarks\label{s4}}

In this paper two approaches to stable differentiation of noisy 
functions are discussed.
The advantage of the first approach is that it contains 
neither matrix inversion nor solving of linear algebraic systems. 
Its computation time is very small. The drawback of the method is that
there is no known {\it a posteriori} choice of $\alpha(\delta)$. 
The second approach is an implementation of the DSM. 
It competes favorably with the VR 
in both computation time and accuracy. 
The DSM competes favorably with the VR 
in solving linear ill-conditioned algebraic systems. 
{\it A posteriori} choice of $\alpha$, an efficient way to compute 
integral \eqref{eq4} for the first method, 
and an efficient discretization of the Volterra equation \eqref{eq1}
with the implementation of the DSM are planned for future research.

\end{document}